\documentclass{amsart}
\usepackage{graphicx}
\usepackage{amssymb}
\usepackage{amsmath}
\newtheorem{theorem}{Theorem}[section]

\newtheorem{corollary}[theorem]{Corollary}

\newtheorem{remark}[theorem]{Remark}

\newenvironment{pot2.1}{{\bf Proof.}}{\hfill\fbox{}\par\vspace{.2cm}}

\numberwithin{equation}{section}

\def\charf {\mbox{{\text 1}\kern-.24em {\text l}}}

\def\bea{\begin{eqnarray*}}
\def\eea{\end{eqnarray*}}
\def\be{\begin{eqnarray}}
\def\ee{\end{eqnarray}}

\begin{document}

\title{Myers-type compactness theorem with the Bakry-Emery Ricci tensor}
\author{Seungsu Hwang}
\address{Department of Mathematics, Chung-Ang University, 84 Heukseok-ro, Dongjak-gu, Seoul 06974, Republic of Korea}
\curraddr{}
\email{seungsu@cau.ac.kr}
\thanks{}

%    author two information
\author{Sanghun Lee}
\address{Department of Mathematics, Chung-Ang University, 84 Heukseok-ro, Dongjak-gu, Seoul 06974, Republic of Korea}
\curraddr{}
\email{kazauye@cau.ac.kr}
% \thanks{Seungsu Hwang was supported by the National Research Foundation of Korea(NRF-2018R1D1A1B05042186).}


\date{}

\dedicatory{}

\begin{abstract}
 In this paper, we first prove the $f$-mean curvature comparison in a smooth metric measure space when the Bakry–Emery Ricci tensor is bounded from below and $|f|$ is bounded. Based on this, we define a Myers–type compactness theorem by generalizing the results of Cheeger, Gromov, and Taylor and of Wan for the Bakry–Emery Ricci tensor. Moreover, we improve a result from Soylu by using a weaker condition on a derivative $f'(t)$.
\end{abstract}

\maketitle
\section{Introduction}
The Myers theorem is a fundamental result in Riemannian geometry; it states that if an $n$-dimensional complete Riemannian manifold $(M,g)$ satisfies $Ric \geq (n-1)H $ with $H>0$, then $M$ is compact and ${\rm diam}(M) \leq \frac{\pi}{\sqrt{H}}$. Here, $Ric$ is the Ricci curvature of the metric $g$. This result has been generalized by different approaches (\cite{AB}, \cite{GL}, \cite{CGM}, and \cite{CS}). Herein, we recall two of them. The first was derived by W. Ambrose.

\begin{theorem}[\cite{AB}] 
 Suppose there exists a point $p$ in an $n$-dimensional complete Riemannian manifold $M$ for which every geodesic $\gamma(t)$ emanating from $p$ satisfies
\[ \int^{\infty}_{0} Ric(\gamma'(s),\gamma'(s)) ds = \infty. \]
Then $M$ is compact.
\end{theorem} 

 The second theorem was derived by J. Cheeger, M. Gromov, and M. Taylor.

\begin{theorem} [\cite{CGM}] \label{thm1.1}
Let $M$ be an $n$-dimensional complete Riemannian manifold. If there exists a fixed point $p \in M$ and $r_{0}, \nu > 0$ such that
$$
Ric(x) \geq (n-1)\frac{(\frac{1}{4} + \nu^{2})}{r^{2}(x)}
$$
holds for all $r(x) \geq r_{0} > 0$, where $r(x)$ is a distance function defined with respect to a fixed point $p \in M$, i.e., $r(x) = d(p,x)$, then $M$ is compact and $\mbox{\rm diam}(M) < r_{0}e^{\frac{\pi}{\nu}}$.
\end{theorem}

 Theorem \ref{thm1.1} was improved by J. Wan as follows:

\begin{theorem} [\cite{JW}] \label{thm1.2}
 Let $M$ be an $n$-dimensional complete Riemannian manifold. If there exists a fixed point $p \in M$, $b \geq 2$, and $r_{0} > 0$ such that
$$
Ric(x) \geq \frac{C(n, b, r_{0})}{(r_{0} + r(x))^{b}} 
$$
for all $r(x) \geq 0$, where $r(x) = d(p,x)$ and $C(n, b, r_{0})$ is a constant depending on $n$, $b$, and $r_{0}$, then $M$ is compact. In particular, $C(n, b, r_{0})$ can be chosen to be equal to $(n-1)\frac{(b-1)^{b}}{(b-2)^{b-2}}r^{b-2}_{0}$ for $b>2$ and $(n-1)(1 + \frac{r_{0}}{\epsilon})$, $\epsilon > 0$ for $b=2$.
\end{theorem}
 
 Now, we consider the $k$-Bakry–Emery Ricci tensor as follows:
\[ Ric_{\,V}^{\,k} := Ric + \frac{1}{2} \mathcal{L}_{V}g - \frac{1}{k} V^{*} \otimes V^{*} \]
for some number $k > 0$, where $\mathcal{L}_{V}$ denotes the Lie derivative in the direction of $V$, and $V^{*}$ is the metric dual of $V$. When $k = \infty$, \[Ric_{\,V}:= Ric_{\,V}^{\,\infty} = Ric + \frac{1}{2}\mathcal{L}_{V}g. \]
In particular, if \[Ric_{\,V} = \lambda g\] for some $\lambda \in \mathbb{R}$, then $(M,g)$ is a Ricci soliton, which is a self-similar solution to the Ricci flow. It is classified as expanding, steady, or shrinking when $\lambda < 0$, $\lambda = 0$, or $\lambda > 0$, respectively. When $V = \nabla f$ for some $f \in C^{\infty}(M)$, the $k$-Bakry--Emery Ricci tensor becomes \[ Ric_{f}^{\,k} := Ric + \mbox{\rm Hess}f  - \frac{1}{k} df \otimes df, \] and the Ricci soliton becomes a gradient Ricci soliton \[Ric_{f} := Ric + {\rm Hess}\, f, \] where ${\rm Hess}\, f$ is the Hessian of $f$.

 G. Wei and W. Wylie have generalized the Myers theorem by using $Ric_{f}$ in \cite{GW}. Several works have attempted to generalize this result (\cite{JY}, \cite{MR}, \cite{HT}, and \cite{SZ}). Note that H. Qiu generalized Theorem \ref{thm1.2} to Theorem \ref{thm1.3} by using $Ric_{\, V}$.

\begin{theorem} [\cite{HQ}] \label{thm1.3}
 Let $M$ be an $n$-dimensional complete Riemannian manifold, $V$ be a smooth vector field on $M$, and $h: [0, \infty) \rightarrow (0, \infty)$ be a continuous function. Let $r(x) = d(p,x)$ be the distance function from $p \in M$. Assume that $\langle V,\nabla r \rangle \leq \delta_{1}$ along a minimal geodesic from each point $\tilde{p} \in M$; here, $\delta_{1}$ is a constant. Suppose
$$
Ric_{\, V}(x) \geq \delta_{2}h(r(x)),
$$ 
where $\delta_{2}$ is a positive constant depending only on $h$, $n$, and $\delta_{1}$, then $M$ is compact. In particular, $\delta_{2}$ can be chosen as $\left(\frac{n-1}{\epsilon} + 2\delta_{1} \right)\cdot \left(\frac{1}{\int^{\infty}_{\epsilon}h(s)ds} \right) + \epsilon_{1}$ ($\epsilon$ and $\epsilon_{1}$ are arbitrary positive constants).
\end{theorem}
\vspace{.2cm}
 
 In this paper, we improve Theorem \ref{thm1.3} by weakening the condition in the case of $Ric_{f}$ (Theorem \ref{thm1} and \ref{thm2}). Here, instead of the condition $|f| \leq k$, our conditions are $|f| \leq \delta(d(p,x) + 1)$ ($\delta > 0$) for Theorem \ref{thm1} and $\partial_{t}f \geq -a$ ($a \geq 0$) for Theorem \ref{thm2}. Additionally, we prove Theorem \ref{thm3} in the case of $Ric_{f}^{\, k}$.

\begin{theorem}\label{thm1}
 Let $M$ be an $n$-dimensional complete Riemannian manifold and $h: [0, \infty) \rightarrow (0,\infty)$ be a continuous function. If there exists $p \in M$ and $|f| \leq \delta (d(p,x) + 1)$ for some constant $\delta > 0$ such that
\be \label{eq1}
Ric_{f}(x) \geq C_{1}(h, n, \delta, \epsilon)h(r(x)),
\ee
where $r(x) = d(p,x)$ and $C_{1}$ is a positive constant depending only on $h$, $n$, $\delta$, and $\epsilon$, then $M$ is compact. In particular, $C_{1}$ can be chosen as \\
$\left(4\delta + \frac{n + 4\delta (\epsilon + 1) - 1}{\epsilon} \right)\cdot \left(\frac{1}{\int^{\infty}_{\epsilon}h(s)ds}\right ) + \epsilon_{1}$ ($\epsilon$ and $\epsilon_{1}$ are arbitrary positive constants).
\end{theorem}

 As a corollary, we make the following improvement to Theorem \ref{thm1.2}.

\begin{corollary}\label{coro1}
 Let $M$ be an $n$-dimensional complete Riemannian manifold. If there exists $p \in M$, $b \in \mathbb{R}$, $r_{0} >0$, and $|f| \leq \delta (d(p,x) + 1)$ for some constant $\delta > 0$ such that
\be \label{eq2}
Ric_{f}(x) \geq \frac{C_{2}(n, b, r_{0}, \delta, \epsilon)}{(r_{0} + r(x))^{b}}
\ee
for all $r(x) \geq 0$, where $r(x) = d(p,x)$ and $C_{2}$ is a constant depending on $n$, $b$, $r_{0}$, $\delta$, and $\epsilon$, then $M$ is compact. In particular, $C_{2}$ can be chosen as $\left(4\delta + \frac{n + 4\delta(\epsilon + 1) - 1}{\epsilon} \right)\cdot \left((b-1)(r_{0} + \epsilon)^{b-1} \right) + \epsilon_{1}$ for $b>1$, and $C_{2}$ can be chosen as $\epsilon_{1}$ for $b \leq 1$. Here, $\epsilon$ and $\epsilon_{1}$ are arbitrary positive constants.
\end{corollary}

 Furthermore, we have

\begin{theorem}\label{thm2}
 Let $M$ be an $n$-dimensional complete Riemannian manifold and $h: [0, \infty) \rightarrow (0, \infty)$ be a continuous function. If there exists $p \in M$ and $\partial_{t}f \geq -a \, (a \geq 0)$ such that
\be \label{eq3}
Ric_{f}(x) \geq C_{3}(h, n, a)h(r(x)),
\ee
where $r(x) = d(p,x)$ and $C_{3}$ is a positive constant depending on $h$, $n$, and $a$, then $M$ is compact.
In particular, $C_{3}$ can be chosen as $\left(a + \frac{n-1}{\epsilon} \right) \cdot \left(\frac{1}{\int^{\infty}_{\epsilon} h(s) ds}\right) + \epsilon_{1}$ ($\epsilon$ and $\epsilon_{1}$ are arbitrary positive constants).
\end{theorem}

 As a corollary, we obtain the following result, similar to Corollary \ref{coro1}.

\begin{corollary}\label{coro2}
 Let $M$ be an $n$-dimensional complete Riemannian manifold. If there exists $p \in M$, $b \in \mathbb{R}$, $r_{0} >0$, and $\partial_{t}f \geq -a \, (a \geq 0)$ such that
\be \label{eq4}
Ric_{f}(x) \geq \frac{C_{4}(n,b,r_{0},a)}{(r_{0} + r(x))^{b}}
\ee
for all $r(x) \geq 0$, where $r(x) = d(p,x)$ and $C_{4}$ is a constant depending only on $n$, $b$, $r_{0}$, and $a$, then $M$ is compact.
In particular, $C_{4}$ can be chosen as $(2ar_{0} + (n-1)(b-2))r^{b-2}_{0}\frac{(b-1)^{b}}{(b-2)^{b-1}}$ for $b>2$, $C_{4}$ can be chosen as $(2a + \frac{n-1}{\epsilon})(b-1)(r_{0} + \epsilon)^{b-1}$ for $1 < b \leq 2$, and $C_{4}$ can be chosen as $\epsilon_{1}$ for $b \leq 1$. Here, $\epsilon$ and $\epsilon_{1}$ are arbitrary positive constants.
\end{corollary}

 For the case of $Ric_{f}^{\, k}$, we have

\begin{theorem} \label{thm3}
 Let $M$ be an $n$-dimensional complete Riemannian manifold and $h: [0, \infty) \rightarrow (0, \infty)$ be a continuous function. If there exists $p \in M$ such that
\be \label{eq5}
Ric_{f}^{\, k}(x) \geq C_{5}(h, n+k)h(r(x)),
\ee
where $k \in (0, \infty)$, $r(x) = d(p,x)$, and $C_{5}$ is a positive constant depending on $h$ and $n+k$, then $M$ is compact. In particular, $C_{5}$ can be chosen as\\
$\left(\frac{n+k-1}{\epsilon} \right)\cdot \left(\frac{1}{\int^{\infty}_{\epsilon} h(s) ds} \right) + \epsilon_{1}$ ($\epsilon$ and $\epsilon_{1}$ are arbitrary positive constants).
\end{theorem}

\begin{corollary}\label{coro3}
 Let $M$ be an $n$-dimensional complete Riemannian manifold. If there exists $p \in M$, $b \in \mathbb{R}$, and $r_{0} > 0$ such that
\be \label{eq6}
Ric_{f}^{\, k}(x) \geq \frac{C_{6}(n+k,b,r_{0})}{(r_{0} + r(x))^{b}} 
\ee
for all $r(x) \geq 0$, where $k \in (0, \infty)$, $r(x) = d(p,x)$, and $C_{6}$ is a constant depending on $n+k$, $b$, and $r_{0}$, then $M$ is compact. In particular, $C_{6}$ can be chosen as $(n+k-1)\frac{(b-1)^{b}}{(b-2)^{b-2}}r_{0}^{b-2}$ for $b>2$, $C_{6}$ can be chosen as\\
 $(\frac{n+k-1}{\epsilon})(b-1)(r_{0} + \epsilon)^{b-1}$ for $1 < b \leq 2$, and $C_{6}$ can be chosen as $\epsilon_{1}$ for $b \leq 1$. Here, $\epsilon$ and $\epsilon_{1}$ are arbitrary positive constants. 
\end{corollary}

\begin{remark}
 If $\int^{\infty}_{\epsilon} h(s) ds = \infty$, then the constants $C_{1}$, $C_{3}$, and $C_{5}$ can be chosen as arbitrary positive real numbers.
\end{remark}

 Moreover, we obtain an Ambrose-type result. Theorem 1.1 was generalized by S. Zhang \cite{SZ} using the Bakry–Emery Ricci tensor. Another generalization was derived by M.P. Cavalcate, J.Q. Oliveira, and M.S. Santos \cite{MP} with the condition $f$: $\frac{df}{dt} \leq 0$. The following is a further improvement made by Y. Soylu \cite{YS}.

\begin{theorem}[\cite{YS}]
 Let $M$ be an $n$-dimensional complete Riemannian manifold, where $n\geq 2$. Suppose there exists a point $p \in M$ such that every geodesic $\gamma(t)$ emanating from $p$ satisfies
\[ \int^{\infty}_{0} Ric_{f}(\gamma'(t), \gamma'(t)) dt = \infty ,\]
and $f'(t) \leq \frac{1}{4}(1 - \frac{1}{t})$ for all $t \geq 1$. Then, $M$ is compact.
\end{theorem}

 In this paper, we improve the above result by finding a weaker condition for the derivative $f'(t)$:

\begin{theorem}\label{thm4}
 Let $M$ be an $n$-dimensional complete Riemannian manifold, where $n \geq 2$. Suppose there exists a point $p \in M$ such that every geodesic $\gamma(t)$ emanating from $p$ satisfies
$$
\int^{\infty}_{0} Ric_{f}(\gamma'(t), \gamma'(t)) dt = \infty ,
$$
and $f'(t) \leq C(1 - \frac{1}{t^{\alpha}})$ for some constants $C > 0$, $\alpha > 1$, and for all  $t \geq 1$. \\ 
Then, $M$ is compact.
\end{theorem}

 This paper is organized as follows. In Section 2, we prove the $f$-mean curvature comparison theorem. In Section 3, we prove Theorems \ref{thm1}, \ref{thm2}, and \ref{thm3}. In Section 4, we prove Corollaries \ref{coro1}, \ref{coro2}, and \ref{coro3}. In the final section, we prove Theorem \ref{thm4}.

\section{$f$-Mean curvature comparison}
 
 We first recall a few definitions. Let $(M,g, e^{-f}dv_{g})$ be a smooth metric measure space on  
 an $n$-dimensional complete Riemannian manifold $M$. For the measure $e^{-f}dv_{g}$, the $f$-mean curvature is defined by $m_{f} := m - \partial_{r}f$, where $m$ is the mean curvature of the geodesic sphere with an inward-pointing normal vector. Then, the $f$-Laplacian is defined by $\Delta_{f} := \Delta - \langle \nabla f, \nabla \rangle$. Note that $m_{f}(r) = \Delta_{f}(r)$ and $m(r) = \Delta (r)$, where $r$ is the distance function. Now, we prove the $f$-mean curvature comparison theorems.

\begin{theorem} \label{thm2.1}
 Let $M$ be an $n$-dimensional complete Riemannian manifold with $Ric_{f}(\gamma',\gamma') \geq (n-1)H$, $H \in \mathbb{R}$. Fix $p \in M$. If $|f| \leq \delta(d(p,x) + 1)$ for some constant $\delta > 0$ and $d(p,x)$ is the distance function from $p$ to $x$, then
\be \label{eq2.1}
m_{f}(t) \leq m^{n + 4\delta(t + 1)}_{H}(t)
\ee
(when $H>0$, assume t $\leq \frac{\pi}{4\sqrt{H}}$). \\
Here, $m^{n + 4\delta(t+1)}_{H}(t)$ is the mean curvature of the geodesic sphere in $M^{n+4\delta(t+1)}_{H}$. The simply connected model space of dimension $n+4\delta(t+1)$ has a constant curvature $H$, and $m_{H}$ is the mean curvature of the model space of dimension $n$.
\end{theorem}
\begin{pot2.1}
By the Bochner formula and Schwarz inequality, the distance function $r$ satisfies the Riccati inequality,
$$
- \frac{(\Delta r)^{2}}{n-1} \geq \frac{\partial}{\partial r}(\Delta r) + Ric(\nabla r, \nabla r).
$$ 
i.e.,
$$
m'(t) \leq - \frac{m^{2}(t)}{n-1} - Ric(\gamma'(t),\gamma'(t)).
$$
The equality holds if and only if the radial sectional curvatures are constant. Therefore, the mean curvature of the model space $m_{H}$ satisfies
$$
m'_{H}(t) = -\frac{m^{2}_{H}(t)}{n-1} - (n-1)H.
$$
Since $Ric_{f}(\gamma',\gamma') = Ric(\gamma',\gamma') +{\rm Hess}\, f(\gamma',\gamma')$, we obtain
\be \label{eq2.2}
(m(t) - m_{H}(t))' \leq -\frac{m^{2}(t) - m^{2}_{H}(t)}{n-1} + \nabla\nabla f(\gamma'(t),\gamma'(t)).
\ee

 Let $sn_{H}(t)$ be the solution to $sn''_{H}(t) + Hsn_{H}(t) = 0$ such that $sn_{H}(0) = 0$ and $sn'_{H}(0) = 1$, then
\be \label{eq2.4}
m_{H}(t) = (n-1)\frac{sn'_{H}(t)}{sn_{H}(t)}.
\ee
We compute
$$
(sn_{H}^{2}(m - m_{H}))' = 2sn_{H}sn'_{H}(m - m_{H}) + sn_{H}^{2}(m - m_{H})'.
$$
By inequality (\ref{eq2.2}), we have
$$
(sn_{H}^{2}(m - m_{H}))' \leq sn_{H}^{2} \left (-\frac{(m - m_{H})^{2}}{n-1} + \nabla\nabla f \right) 
\le sn_{H}^{2}\nabla\nabla f.
$$
Since $\nabla\nabla f(\gamma'(t),\gamma'(t)) = \frac{\partial}{\partial t}\langle \nabla f, \gamma'\rangle(t)$, integrating the above inequality from $0$ to $t$ yields
$$
sn_{H}^{2}(t)m(t) \leq sn_{H}^{2}(t)m_{H}(t) + \int^{t}_{0} sn_{H}^{2}(s)\frac{\partial}{\partial s}
\langle \nabla f, \gamma\rangle(s) ds.$$
By performing integration by parts on the last term, we obtain
\bea
sn_{H}^{2}(t)m(t) & \leq & sn_{H}^{2}(t)m_{H}(t) + sn_{H}^{2}(t)\langle \nabla f, \gamma'(t)\rangle \\
& &- \int^{t}_{0}(sn_{H}^{2}(s))'\langle \nabla f, \gamma'\rangle(s) ds.
\eea
Since $\langle \nabla f, \gamma'\rangle (t) = \frac{\partial}{\partial t}(f(\gamma(t)))$ and $m_{f}(t) = m(t) - \frac{\partial}{\partial t}f(\gamma(t))$, we have
$$
sn_{H}^{2}(t)m_{f}(t) \leq sn_{H}^{2}(t)m_{H}(t) - \int^{t}_{0}(sn_{H}^{2}(s))'\frac{\partial}{\partial s}(f(\gamma(s))) ds.
$$
Conducting integration by parts on the last term again, we obtain
$$
sn_{H}^{2}(t)m_{f}(t) \leq sn_{H}^{2}(t)m_{H}(t) - f(\gamma(t))(sn_{H}^{2}(t))' + \int^{t}_{0} (sn_{H}^{2}(s))''f(\gamma(s))ds.$$
If $|f| \leq \delta(d(p,x) + 1)$ and $t \in (0, \frac{\pi}{4\sqrt{H}}]$ when $H > 0$, then $(sn_{H}^{2}(t))' \geq 0$ and $(sn_{H}^{2}(t))'' \geq 0$. Therefore, we have
\bea
sn_{H}^{2}(t)m_{f}(t) &\leq& sn_{H}^{2}(t)m_{H}(t) + \delta(t+1)(sn_{H}^{2}(t))' + \delta\int^{t}_{0}(sn_{H}^{2}(s))''(s+1) ds \\
&=& sn_{H}^{2}(t)m_{H}(t) + \delta(t+1)(sn_{H}^{2}(t))' + \delta\int^{t}_{0}(sn_{H}^{2}(s))''s \, ds \\
& & + \delta\int^{t}_{0}(sn_{H}^{2}(s)'' ds \\
&=& sn_{H}^{2}(t)m_{H}(t) + 2\delta(t+1)(sn_{H}^{2}(t))' - \delta sn_{H}^{2}(t).
\eea
From (\ref{eq2.4}), we compute $$(sn_{H}^{2}(t))' = 2sn_{H}'(t)sn_{H}(t) = \frac{2sn_{H}^{2}(t)m_{H}(t)}{n-1}.$$ Thus, we can see that
$$
sn_{H}^{2}(t)m_{f}(t)  \leq  sn_{H}^{2}(t)m_{H}(t) + \frac{4\delta(t+1)sn_{H}^{2}(t)m_{H}(t)}{n-1} - \delta sn_{H}^{2}(t).$$
Therefore,
\bea
m_{f}(t) &\leq& m_{H}(t) + \frac{4\delta(t+1)m_{H}(t)}{n-1} - \delta \\
&\leq& m_{H}(t)\left(1 + \frac{4\delta(t+1)}{n-1}\right) = m_{H}^{n + 4\delta(t+1)}(t).
\eea
\end{pot2.1}

\begin{theorem} [\cite{GW}] \label{thm2.2}
 Let $M$ be an $n$-dimensional complete Riemannian manifold with $Ric_{f}(\gamma',\gamma') \geq (n-1)H$, $H \in \mathbb{R}$. Fix $p \in M$. If $\partial_{t}f \geq -a$ for some constant $a \geq 0$, then
\be \label{eq2.3}
m_{f}(t) \leq m_{H}(t) + a
\ee
(when $H > 0$ assume $t \leq \frac{\pi}{2\sqrt{H}}$).
\end{theorem}

 The proof of Theorem \ref{thm2.2} is similar to the above proof. For the detailed proof, see Theorem 1.1 of \cite{GW}.

\section{Proof of Theorems \ref{thm1}, \ref{thm2}, and \ref{thm3}}

 In this section, we prove Theorem \ref{thm1}, \ref{thm2}, and \ref{thm3}. The proof uses the $f$-mean curvature comparison calculated above and Riccati inequality. \\

 First, we prove Theorem \ref{thm1}. Suppose $M$ is non-compact. For any $p \in M$, there exists a unit speed ray $\gamma(t)$ starting from $p$. Let $r(x) = d(p,x)$ be the distance function from $p$. By the Bochner formula and Schwarz inequality, we have
$$
\frac{m^{2}(t)}{n-1} \leq -m_{f}'(t) - Ric_{f}(\gamma', \gamma').
$$
Integrating the above inequality from $\epsilon$ to $t$ we obtain
$$
\frac{1}{n-1}\int^{t}_{\epsilon}m^{2}(s) ds \leq - \int^{t}_{\epsilon}m_{f}'(s) ds - \int^{t}_{\epsilon} Ric_{f}(\gamma',\gamma') ds. 
$$
From (\ref{eq1}), the above inequality becomes
\be \label{eq3.1}
\frac{1}{n-1}\int^{t}_{\epsilon}m^{2}(s) ds \leq -m_{f}(t) + m_{f}(\epsilon) - C_{1}\int^{t}_{\epsilon} h(s) ds.
\ee
Now, we claim that $-4\delta \leq m_{f}(t) \leq \frac{n+4\delta(t+1)-1}{t}$. \\
Since $Ric_{f} > 0$, $m_{f}(t) \leq \frac{n+4\delta(t+1)-1}{t}$ holds by inequality (\ref{eq2.1}), it suffices to show that $-4\delta \leq m_{f}(t)$. 

 Consider the excess function
 $$
e(x) = d(p,x) + d(x,\gamma(i)) - i
$$
for $0 \leq t \leq i$. By the triangle inequality, we have
$$
e(x) \geq 0 \,\,\, {\rm and} \,\,\, e(\gamma(t)) = 0.
$$
Therefore,
$$\Delta_{f}(e)(\gamma(t)) \geq 0.
$$
It follows that
$$
\Delta_{f}(e)(\gamma(t)) = \Delta_{f}d(p, \gamma(t)) + \Delta_{f}d(\gamma(t), \gamma(i)) \geq 0.
$$
i.e.,
$$
\Delta_{f}d(p, \gamma(t)) = m_{f}(t) \geq -\Delta_{f}d(\gamma(t),\gamma(i)) = -m_{f}(i-t).
$$
Thus,
$$
m_{f}(t) \geq -m_{f}(i-t) \geq -\frac{n+ 4\delta((i-t) + 1) -1}{i-t}.
$$
If $i \rightarrow \infty$, then $m_{f}(t) \geq -4\delta$. Hence, we proved the claim that
\be \label{eq3.2}
-4\delta \leq m_{f}(t) \leq \frac{n+4\delta(t+1)-1}{t}.
\ee
From (\ref{eq3.1}) and (\ref{eq3.2}), we derive
$$
\frac{1}{n-1}\int^{t}_{\epsilon}m^{2}(s) ds \leq 4\delta + \frac{n+ 4\delta(\epsilon + 1) -1}{\epsilon} - C_{1}\int^{t}_{\epsilon} h(s) ds.
$$
Let $t \rightarrow \infty$. Then, the above inequality becomes
\be \label{eq3.3}
0 \leq \frac{1}{n-1}\int^{\infty}_{\epsilon}m^{2}(s) ds \leq 4\delta + \frac{n+ 4\delta(\epsilon + 1) -1}{\epsilon} - C_{1}\int^{\infty}_{\epsilon} h(s) ds.
\ee
Therefore, we obtain
$$
C_{1}\int^{\infty}_{\epsilon} h(s) ds \leq 4\delta + \frac{n+ 4\delta(\epsilon + 1) -1}{\epsilon}.
$$
Choosing
$$
C_{1} = \left(4\delta + \frac{n + 4\delta (\epsilon + 1) - 1}{\epsilon}\right )\cdot \left(\frac{1}{\int^{\infty}_{\epsilon}h(s)ds}\right ) + \epsilon_{1},
$$
where $\epsilon_{1}>0$. This yields a contradiction. Thus, $M$ must be compact.

\hfill\fbox{}\par\vspace{.5cm}
 
 Second, we prove Theorem \ref{thm2}. Its proof is similar to the previous proof; thus, the setting is identical. Suppose that $M$ is non-compact. Based on the previous proof, we know that
$$
\frac{1}{n-1}\int^{t}_{\epsilon}m^{2}(s) ds \leq - \int^{t}_{\epsilon}m_{f}'(s) ds - \int^{t}_{\epsilon} Ric_{f}(\gamma',\gamma') ds.
$$
From (\ref{eq3}), we have
$$ 
\frac{1}{n-1}\int^{t}_{\epsilon}m^{2}(s) ds \leq -m_{f}(t) + m_{f}(\epsilon) - C_{3}\int^{t}_{\epsilon} h(s) ds.
$$
We claim that $-a \leq m_{f}(t) \leq \frac{n-1}{t} + a$. Based on Theorem \ref{thm2.2}, $m_{f}(t) \leq \frac {n-1}{t} + a$ and $-a \leq m_{f}(t)$ hold according to the same argument as in the previous proof. Thus,
$$
\frac{1}{n-1}\int^{t}_{\epsilon} m^{2}(s) ds \leq 2a + \frac{n-1}{\epsilon} - C_{3}\int^{t}_{\epsilon} h(s) ds.
$$
Let $t \rightarrow \infty$. Then we obtain
$$
0 \leq \frac{1}{n-1}\int^{\infty}_{\epsilon} m^{2}(s) ds \leq 2a + \frac{n-1}{\epsilon} - C_{3}\int^{\infty}_{\epsilon} h(s) ds.
$$
If we choose
$$
C_{3} = \left(2a + \frac{n-1}{\epsilon}\right)\cdot \left(\frac{1}{\int^{\infty}_{\epsilon} h(s) ds} \right) + \epsilon_{1}
$$
for some positive constant $\epsilon_{1}$, then this yields a contradiction. Thus, $M$ must be compact. 
\hfill\fbox{}\par\vspace{.5cm}
 
 Finally, we prove Theorem \ref{thm3}. The setting is identical to that of the previous proof. Suppose that $M$ is non-compact. By the Bochner formula and Schwarz inequality for the distance function $r$, we have the following Riccati inequality (Appendix A of \cite{GW}):
\be \label{eq3.4}
\frac{m_{f}^{2}(t)}{k + n -1} \leq -m_{f}'(t) - Ric_{f}^{\, k}(\gamma'(t), \gamma'(t)).
\ee
Integrating both sides of (\ref{eq3.4}), we obtain
\bea
\frac{1}{k + n -1} \int^{t}_{\epsilon} m_{f}^{2}(s) ds \leq -m_{f}(t) + m_{f}(\epsilon) - \int^{t}_{\epsilon} Ric_{f}^{\, k}(\gamma',\gamma') ds.
\eea
From (\ref{eq5}), the above inequality becomes
\bea
\frac{1}{k + n -1} \int^{t}_{\epsilon} m_{f}^{2}(s) ds \leq -m_{f}(t) + m_{f}(\epsilon) - C_{5}\int^{t}_{\epsilon}h(s) ds.
\eea
We claim that $0 \leq m_{f}(t) \leq \frac{n + k -1}{t}$. By Theorem A.1 in \cite{GW}, $m_{f}(t) \leq \frac{n + k - 1}{t}$ and $m_{f}(t) \geq 0$ hold by the same argument as in the previous proof. Thus
\bea
\frac{1}{k + n -1} \int^{t}_{\epsilon} m_{f}^{2}(s) ds \leq \frac{n+k-1}{\epsilon} - C_{5}\int^{t}_{\epsilon} h(s) ds.
\eea
Let $t \rightarrow \infty$. Then, we have
\bea
0 \leq \frac{1}{k + n -1} \int^{\infty}_{\epsilon} m_{f}^{2}(s) ds \leq \frac{n+k-1}{\epsilon} - C_{5}\int^{\infty}_{\epsilon} h(s) ds.
\eea
Choosing
\bea
C_{5} = \left(\frac{n+k-1}{\epsilon} \right) \cdot \left(\frac{1}{\int^{\infty}_{\epsilon}h(s) ds} \right) + \epsilon_{1},
\eea
where $\epsilon_{1} > 0$. This is a contradiction. Thus, $M$ is compact.
\hfill\fbox{}\par\vspace{.5cm}

\section{Proof of Corollaries \ref{coro1}, \ref{coro2}, and \ref{coro3}}

 In this section, we prove Corollaries \ref{coro1}, \ref{coro2}, and \ref{coro3}. These corollaries generalize Theorem \ref{thm1.2} by using the Bakry–Emery Ricci tensor. The settings of the corollaries are same as those of the previous theorems. \\ 
 
 First, we prove Corollary \ref{coro1}. Suppose that $M$ is non-compact. We know that
$$
\frac{1}{n-1}\int^{t}_{\epsilon}m^{2}(s) ds \leq - \int^{t}_{\epsilon}m_{f}'(s) ds - \int^{t}_{\epsilon} Ric_{f}(\gamma',\gamma') ds.
$$
By assumption (\ref{eq2}) and $-4\delta \leq m_{f}(t) \leq \frac{n + 4\delta(t+1) -1}{t}$, the above inequality becomes
\be \label{eq4.1}
\frac{1}{n-1}\int^{t}_{\epsilon}m^{2}(s) ds &\leq& 4\delta + \frac{n + 4\delta(\epsilon + 1)-1}{\epsilon} - C_{2} \int^{t}_{\epsilon} \frac{1}{(r_{0} + s)^{b}} ds.
\ee
Now, we obtain $C_{2}$. First, when $b>1$, we have
$$
\frac{1}{n-1}\int^{t}_{\epsilon}m^{2}(s) ds \leq 4\delta + \frac{n + 4\delta(\epsilon + 1)-1}{\epsilon} - \frac{C_{2}}{b-1} \left(\frac{1}{(r_{0} + \epsilon)^{b-1}} - \frac{1}{(r_{0} + t)^{b-1}} \right).
$$
If $t \rightarrow \infty$, then we obtain
$$
0 \leq \frac{1}{n-1}\int^{\infty}_{\epsilon}m^{2}(s) ds \leq 4\delta + \frac{n + 4\delta(\epsilon + 1)-1}{\epsilon} - \frac{C_{2}}{b-1} \left(\frac{1}{(r_{0} + \epsilon)^{b-1}} \right).
$$
Choosing 
$$C_{2} = \left(4\delta + \frac{n + 4\delta(\epsilon + 1) - 1}{\epsilon} \right)\cdot \left((b-1)(r_{0} + \epsilon)^{b-1} \right) + \epsilon_{1}$$ for any positive constant $\epsilon_{1}$ yields a contradiction. \\
Second, when $b=1$, inequality (\ref{eq4.1}) becomes
$$
\frac{1}{n-1}\int^{t}_{\epsilon}m^{2}(s) ds \leq 4\delta + \frac{n + 4\delta(\epsilon + 1)-1}{\epsilon} - C_{2}\cdot {\rm ln}\left( \frac{r_{0} + t}{r_{0} + \epsilon} \right).
$$
If $t \rightarrow \infty$ and $C_{2} = \epsilon_{1}$, the above inequality is a contradiction. \\
Similarly, when $b<1$, if $t \rightarrow \infty$ and $C_{2} = \epsilon_{1}$ in inequality (\ref{eq4.1}), we obtain a contradiction. Thus, $M$ must be compact.
\hfill\fbox{}\par\vspace{.5cm}
 
 Secondly, we prove Corollary \ref{coro2}. Suppose that $M$ is non-compact. By assumption (\ref{eq4}) and $-a \leq m_{f}(t) \leq \frac{n-1}{t} + a$, we have
\be \label{eq4.2}
\frac{1}{n-1}\int^{t}_{\epsilon}m^{2}(s) ds \leq 2a + \frac{n-1}{\epsilon} - C_{4} \int^{t}_{\epsilon} \frac{1}{(r_{0} + s)^{b}} ds.
\ee
Now, we obtain $C_{4}$. First, when $b>2$, we have
$$
\frac{1}{n-1}\int^{t}_{\epsilon}m^{2}(s) ds \leq 2a + \frac{n-1}{\epsilon} - \frac{C_{4}}{b-1}\left(\frac{1}{(r_{0} + \epsilon)^{b-1}} - \frac{1}{(r_{0} + t)^{b-1}} \right).
$$
If $t \rightarrow \infty$, then we obtain
$$
0 \leq \frac{1}{n-1}\int^{\infty}_{\epsilon}m^{2}(s) ds \leq 2a + \frac{n-1}{\epsilon} - \frac{C_{4}}{b-1}\left(\frac{1}{(r_{0} + \epsilon)^{b-1}} \right).
$$
Solving
$$
2a + \frac{n-1}{\epsilon} - \frac{C_{4}}{(b-1)(r_{0} + \epsilon)^{b-1}} \leq 0,
$$
it follows that
\be \label{eq4.3}
C_{4} \geq \left(2a + \frac{n-1}{\epsilon} \right)(b-1)(r_{0} + \epsilon)^{b-1}.
\ee
When $\epsilon = \frac{r_{0}}{b-2}$, $\frac{(r_{0} + \epsilon)^{b-1}}{\epsilon}$ has a minimum value. Substituting this into (\ref{eq4.3}), we obtain
$$
C_{4} \geq (2ar_{0} + (n-1)(b-2))r_{0}^{b-2}\frac{(b-1)^{b}}{(b-2)^{b-1}}.
$$
Then, we can choose $C_{4}(n,b,r_{0},a) = (2ar_{0} + (n-1)(b-2))r_{0}^{b-2}\frac{(b-1)^{b}}{(b-2)^{b-1}}$. Similarly, from inequality (\ref{eq4.3}), we can choose $C_{4} = \left(2a + \frac{n-1}{\epsilon} \right)(b-1)(r_{0} + \epsilon)^{b-1}$ for $1 < b \leq 2$. When $b \leq 1$, it is the same as Corollary \ref{coro1}. Thus, $M$ is compact.
\hfill\fbox{}\par\vspace{.5cm}

 Finally, we prove Corollary \ref{coro3}. Suppose that $M$ is a non-compact. By assumption (\ref{eq6}) and $0 \leq m_{f}(t) \leq \frac{n + k -1}{t}$, we obtain
\be \label{eq4.4}
\frac{1}{k + n -1} \int^{t}_{\epsilon} m_{f}^{2}(s) ds \leq \frac{n+k-1}{\epsilon} - C_{6}\int^{t}_{\epsilon} \frac{1}{(r_{0} + s)^{b}} ds. 
\ee
Now, we obtain $C_{6}$. First, for $b>2$, we have
\bea
\frac{1}{k + n -1} \int^{t}_{\epsilon} m_{f}^{2}(s) ds \leq \frac{n+k-1}{\epsilon} - \frac{C_{6}}{b-1}\left(\frac{1}{(r_{0} + \epsilon)^{b-1}} - \frac{1}{(r_{0} + t)^{b-1}} \right).
\eea
Let $t \rightarrow \infty$. We obtain
\bea
0 \leq \frac{1}{k + n -1} \int^{t}_{\epsilon} m_{f}^{2}(s) ds \leq \frac{n+k-1}{\epsilon} - \frac{C_{6}}{b-1}\left(\frac{1}{(r_{0} + \epsilon)^{b-1}} \right).
\eea
If $\frac{n+k-1}{\epsilon} - \frac{C_{6}}{b-1}\left(\frac{1}{(r_{0} + \epsilon)^{b-1}} \right) \leq 0$, then we have
\be \label{eq4.5}
C_{6} \geq \frac{n+k-1}{\epsilon}(b-1)(r_{0} + \epsilon)^{b-1}.
\ee
Similar to the above proof, when $\epsilon = \frac{r_{0}}{b-2}$, we can choose $C_{6}(n+k,b,r_{0}) = (n+k-1)\frac{(b-1)^{b}}{(b-2)^{b-2}}r_{0}^{b-2}$. Similarly, from inequality (\ref{eq4.4}) and (\ref{eq4.5}), we can choose $C_{6} = (\frac{n+k-1}{\epsilon})(b-1)(r_{0} + \epsilon)^{b-1}$ for $1< b \leq 2$ and $C_{6} = \epsilon_{1}$ for $b \leq 1$, where $\epsilon_{1}$ is an arbitrary positive constant. Thus, $M$ is compact.
\hfill\fbox{}\par\vspace{.5cm}

\section{Ambrose-type Theorem}
 
 In this section, we prove Theorem~\ref{thm4}. Suppose that $M$ is non-compact.
Let $\gamma(t)$ be a unit speed ray starting from $p$. For every $ t > 0$, let $m(t)$ be the Laplacian of the distance function from a fixed point $p \in M$.
 By applying the Bochner formula
\bea
\frac{1}{2} \Delta |\nabla u|^{2} = |{\rm Hess}\, u|^{2} + \langle \nabla u, \nabla (\Delta u) \rangle + Ric(\nabla u, \nabla u)
\eea
to the distance function $r(x) = d(x,p)$, we obtain
\bea
0 = |{\rm Hess}\, r|^{2} + \frac{\partial}{\partial r}(\Delta r) + Ric(\nabla r, \nabla r).
\eea
By the Schwarz inequality, we have the Riccati inequality
$$
m'(t) + \frac{m^{2}(t)}{n-1} + Ric(\gamma',\gamma') \leq 0.
$$
By adding ${\rm Hess}\, f(\gamma',\gamma')$ to both sides of this inequality, we obtain
\be
Ric_{f}(\gamma',\gamma') \leq -m'(t) - \frac{m^{2}(t)}{n-1} + \frac{\partial}{\partial t} \langle \nabla f, \gamma' \rangle (t). \label{eq 5.1}
\ee
Integrating both sides of (\ref{eq 5.1}), we obtain 
\bea
\int^{t}_{1} Ric_{f}(\gamma'(s),\gamma'(s)) ds &\leq& \int^{t}_{1} -m'(s) ds - \frac{1}{n-1}\int^{t}_{1} m^{2}(s) ds \nonumber \\
& &+ \langle \nabla f, \gamma' \rangle (t) - \langle \nabla f, \gamma' \rangle (1).
\eea
From assumption
\bea
\lim_{t \rightarrow \infty} \int^{t}_{0} Ric_{f}(\gamma'(t),\gamma'(t)) dt = \infty,
\eea
we have
\be
\lim_{t \rightarrow \infty}\left(- m(t) - \frac{1}{n-1} \int^{t}_{1} m^{2}(s) ds + f'(t)\right) = \infty. \label{eq 5.2}
\ee
By multiplying both sides of (\ref{eq 5.2}) with $\frac{1}{n}$, we obtain
\bea
\lim_{t \rightarrow \infty} \left(- \frac{1}{n}m(t) - \frac{1}{n(n-1)} \int^{t}_{1} m^{2}(s) ds + \frac{1}{n} f'(t)\right) = \infty.
\eea
Thus, for any constant $A \geq C + 2$, there exists $t_{1} > 1$ such that for any $t \geq t_{1}$, we have 
\bea
- \frac{1}{n} m(t) - \frac{1}{n(n-1)}\int^{t}_{1}m^{2}(s)ds + \frac{1}{n} f'(t) \geq A.
\eea
By assumption $f'(t) \leq C(1 - \frac{1}{t^{\alpha}})$, we have
\bea
C + 2 &\leq& - \frac{1}{n} m(t) - \frac{1}{n(n-1)}\int^{t}_{1}m^{2}(s)ds + \frac{C}{n} - \frac{C}{nt^{\alpha}} \\
&\leq& - \frac{1}{n} m(t) - \frac{1}{n(n-1)}\int^{t}_{1}m^{2}(s)ds + C.
\eea
In other words,
\be
- \frac{1}{n} m(t) \geq \frac{1}{n(n-1)}\int^{t}_{1}m^{2}(s)ds  + 2. \label{eq 5.3}
\ee
By multiplying both sides of (\ref{eq 5.3}) by $n$, we obtain 
\be
- m(t) \geq \frac{1}{n-1}\int^{t}_{1}m^{2}(s)ds + 2n, \,\, {\rm for \,\, all }\ t \geq t_{1}. \label{eq 10}
\ee
 Now, we consider the increasing sequence $\{ t_{\ell} \}$ defined by
\bea
t_{\ell + 1} = t_{\ell} + 2^{1 - \ell}, \,\, {\rm for }\  \ell \geq 1.
\eea
Note that $\{ t_{\ell} \}$ converges to $T := t_{1} + 2$ as $\ell \rightarrow \infty$.

 We claim that $- m(t) \geq 2^{\ell}n$ for all $t \geq t_{\ell}$.  We prove this claim through induction. If $\ell = 1$, the claim is trivially true from the inequality in (\ref{eq 10}). Now, for all $t \geq t_{\ell + 1}$, by the inequality in (\ref{eq 10}), we have
\bea
- m(t) &\geq& \frac{1}{n-1}\int^{t}_{1} m^{2}(s) ds + 2n  \nonumber \\
&\geq& \frac{1}{n}\int^{t_{\ell +1}}_{t_{\ell}} m^{2}(s) ds  \nonumber \\
&\geq& \frac{1}{n} \cdot 2^{2 \ell}n^{2}(t_{\ell + 1} - t_{\ell}) = 2^{\ell + 1}n,
\eea
proving the claim. Therefore, 
\bea
\lim_{\ell \rightarrow \infty} - m(t_{\ell}) = - m(T) \geq \lim_{\ell \rightarrow \infty} 2^{\ell}n.
\eea
This contradicts the smoothness of $m(t)$, which completes the proof of Theorem~\ref{thm4}. 

\hfill\fbox{}\par

\end{document}